\documentclass[11pt, a4paper]{article}
\usepackage{amssymb}

\usepackage{amsmath}
\usepackage{color}
\usepackage{natbib}

\ifx\pdfoutput\undefined
  \usepackage[dvips]{graphicx}
\else
  \usepackage[pdftex]{graphicx}
\else
\else
\else
\else
\else
\else
\else
\else
\else

\else
\fi 
\begin{document}

\author{\textsc{Marco}~\textsc{Avarucci} \\
SEFeMEQ, University of Rome ``Tor Vergata'' \and \textsc{Domenico} ~\textsc{%
Marinucci} \\
Department of Mathematics, University of Rome ``Tor Vergata''}
\title{\textsc{Polynomial Cointegration among Stationary Processes with Long
Memory\thanks{%
This paper is part of the first author PhD dissertation. We are grateful to
Franco Peracchi for useful comments on an earlier version. }}}
\maketitle

\begin{abstract}
In this paper we consider polynomial cointegrating relationships among
stationary processes with long range dependence. We express the regression
functions in terms of Hermite polynomials and we consider a form of spectral
regression around frequency zero. For these estimates, we establish
consistency by means of a more general result on continuously averaged
estimates of the spectral density matrix at frequency zero.\newline


\end{abstract}

\section{Introduction}

The extension of the standard cointegration paradigm to more general,
fractional circumstances has drawn increasing attention in the time series
literature over the last decade. The possibility of fractional cointegration
was already mentioned in the seminal paper by \cite{eg87}. \cite{rob94} was
the first to establish consistency for narrow-band estimates of fractional
cointegrating relationships in the stationary case. The properties of this
estimator (which has become known as NBLS) were then investigated under
nonstationary circumstances by \cite{mr01}, \cite{rm01,rm03}. \cite%
{chen_hurv03a,chen_hurv03b} considered principal components methods in the
frequency domain, whereas \cite{vel03}, \cite{rob_hualde03} advocate
pseudo-maximum likelihood methods which improve the efficiency of the
estimates and yield standard asymptotic properties. Cointegration among
stationary processes has also been considered, for instance by \cite{mar00}, %
\cite{cn04}. Many other insightful papers on fractional cointegration have
appeared in the literature, for instance \cite{dol_marm04}, \cite{dav02}.

All these papers have focused on the case of linear cointegration.
Nevertheless, the possibility of polynomial cointegrating relationships
seems of practical interest, for instance (but not exclusively) for
applications to financial data. Nonlinear cointegration has been considered
in the literature (most recently by \cite{kmt05}), but only in
non-fractional circumstances, to the best of our knowledge. In this paper,
we shall focus on nonlinear cointegrating relationships among stationary
long memory processes; the restriction to a stationarity framework is made
necessary by the need to exploit the rich machinery of expansions into
Hermite polynomials, an extremely powerful tool to investigate nonlinear
transformations (see for instance \cite{gir_surg_85}, \cite{arcones94}, \cite%
{surgalis03}). Our general setting can be explained as follows. Let $%
\{A_{t}\}=\{x_{t},e_{t}\},$ $t\in \mathbb{Z}$ be a stationary bivariate time
series with mean zero and covariance such that%
\begin{equation*}
\mathbb{E}A_{t}A_{t+\tau }^{\prime }:=\Gamma (\tau )=\int_{0}^{2\pi
}f(\lambda )e^{i\tau \lambda }d\lambda \text{ ,}
\end{equation*}%
where 
\begin{equation*}
f(\lambda )=\left[ 
\begin{array}{cc}
f_{xx}(\lambda ) & f_{xe}(\lambda ) \\ 
f_{ex}(\lambda ) & f_{ee}(\lambda )%
\end{array}%
\right] \text{ ,}
\end{equation*}%
is the spectral density matrix of $\{A_{t}\}$. We shall take $%
\{x_{t},e_{t}\} $ to be long memory, in the sense that 
\begin{equation}
\gamma _{ab}(\tau )\simeq G_{ab}\tau ^{d_{a}+d_{b}-1}  \label{eq:covariance}
\end{equation}%
for $a,b=x,e$ , $0<d_{a},d_{b}<\frac{1}{2},$ $G_{xx},G_{ee}>0,$ $\left|
G_{xe}\right| \geq 0$. We write $z\sim I(d_{z})$ for long memory processes
with memory parameter $d_{z},\ $and $\simeq $ to denote that the ratio of
the left- and right-hand sides tends to 1.

Now assume there is a polynomial function $g(\cdot )$ such that $\mathbb{E}%
[g(x_{t})]=0$ and 
\begin{equation}
y_{t}=g(x_{t})+e_{t},\qquad 0<d_{e}<d_{y}\leq d_{x}<1/2\text{ ;}
\label{eq:model}
\end{equation}%
in this case, we say that $y_{t},x_{t}$ are nonlinearly cointegrated.
Clearly, the standard (stationary) fractional cointegrating relationship is
obtained in the special case where $g(\cdot )$ is a linear function.

Our main idea in this paper is to write $g(\cdot )$ as a sum of Hermite
polynomials; the coefficients of these polynomials will be estimated by
means of a spectral regression method, as \cite{rob94}, \cite{mar00} and %
\cite{mr01}. We shall show that, by using a degenerating band of frequencies
around the origin, then the estimator of these coefficients is consistent,
even if $x_{t}$ and $e_{t}$ are allowed to be correlated. The plan of this
paper is as follows; in Section 2 we review some results on long memory
processes, Hermite polynomials and the diagram formula; in Section 3 we
discuss consistent estimation of nonlinear cointegrating relationships,
whereas in Section 4 we collect some comments and directions for future
research. Some technical results are collected in an Appendix. In the
sequel, $C$ denotes a generic, positive, finite constant, which need not to
be the same all the time it is used; for two generic matrices $A$ and $B$,
of equal dimension, we say that $A\simeq B$ if, for each $(i,j)$, the ratio
of the $(i,j)$-th elements of $A$ and $B$ tends to unity.

\section{Nonlinear transformation of long memory process}

It is well-known that, under regularity conditions, a consistent estimator
of the spectral density matrix at frequency zero is given by (see for
instance \cite{hannan70}, p.246) 
\begin{equation}
\hat{f}_{ab}(0)=\int_{-\pi }^{\pi }K_{M}(\lambda )I_{ab}(\lambda )d\lambda
\label{eq:one}
\end{equation}%
where the kernel $K_{M}(\lambda )$ is symmetric and such that $\int_{-\pi
}^{\pi }K_{M}(\lambda )d\lambda =1$, and $M$ is a positive integer
satisfying the bandwidth condition $1/M+M/n\rightarrow 0.$ Here, $%
I_{ab}(\lambda )$ is the periodogram, that is 
\begin{equation*}
I_{ab}(\lambda )=\frac{1}{2\pi }\sum_{\tau =-n+1}^{n-1}c_{ab}(\tau )\exp
(-i\lambda \tau )\text{ ,}
\end{equation*}%
for 
\begin{equation*}
c_{ab}(\tau )=\left\{ 
\begin{array}{ll}
n^{-1}\sum_{t=1}^{n-\tau }a_{t}b_{t+\tau } & \qquad \tau \geq 0 \\ 
n^{-1}\sum_{t=|\tau |+1}^{n}a_{t}b_{t-|\tau |} & \qquad \tau <0%
\end{array}%
\right. ,\text{ }a_{t},b_{t}\in \mathbb{R}\text{ , }t=1,2...,n\text{ .}
\end{equation*}%
Equation (\ref{eq:one}) can be rewritten as 
\begin{equation}
\hat{f}_{ab}(0)=\frac{1}{2\pi }\sum_{-n+1}^{n-1}k_{M}(\tau )c_{ab}(\tau )%
\text{ , }k_{M}(\tau ):=\int_{\pi }^{\pi }K_{M}(\lambda )\exp ({i\tau
\lambda )}d\lambda \text{ }.  \label{eq:lagw}
\end{equation}%
As usual, we call the function $k_{M}(\cdot )$ a lag window, and the
corresponding estimator (\ref{eq:lagw}) the lag window spectral density
estimator. The asymptotic behaviour of $\hat{f}_{ab}(0)$ under short range
dependence conditions is now standard textbook material, see for instance %
\cite{tani_kaki00}. Under long memory circumstances, \cite{rob94} and \cite%
{lobato97} investigated the behaviour of, respectively, a discrete
univariate and multivariate version of (\ref{eq:one}); \cite{rob94} propose
an application of this statistic for the estimation of the cointegrating
vector in a stationary framework. The behaviour of the lag window estimator
under long memory is discussed by \cite{mar00}.

In the latter reference, the following linear cointegrating relationship is
considered: 
\begin{equation*}
y_{t}=\beta x_{t}+e_{t},\quad 0\leq d_{e}<d_{y}=d_{x}<1/2\text{ ,}
\end{equation*}%
and it is shown that $\beta $ is consistently estimated by%
\begin{equation*}
\widetilde{\beta }_{M}=\frac{\int_{-\pi }^{\pi }K_{M}(\lambda
)I_{ab}(\lambda )d\lambda }{\int_{-\pi }^{\pi }K_{M}(\lambda )I_{bb}(\lambda
)d\lambda }=\frac{\sum_{\tau =-n+1}^{n-1}k_{M}(\tau )c_{ab}(\tau )}{%
\sum_{\tau =-n+1}^{n-1}k_{M}(\tau )c_{bb}(\tau )}\text{ ,}
\end{equation*}%
under the bandwidth condition $M^{2}=o(n)$ as $n\rightarrow \infty $. We
consider here the same procedure but under generalized, polynomial
circumstances. The investigation of nonlinear transformations requires the
computation of the cumulants of the Hermite polynomials; this is usually
achieved by means of the so-called diagram formula (see for instance \cite%
{arcones94}). We introduce now the main ideas behind this approach. Consider
the polynomial transformation of the Gaussian process (see also \cite{dg02}) 
\begin{equation}
g(z_{t})=\sum_{k=k_{0}}^{K}a_{k}z_{t}^{k}=\sum_{k=k_{0}}^{K}b_{k}H_{k}(z_{t})%
\text{ , }k_{0}\geq 1\text{ , }t=1,2,\dots n\text{ ,}  \label{eq:hermite}
\end{equation}%
where 
\begin{equation*}
b_{k}=\frac{\mathbb{E}(g(z)H_{k}(z))}{k!\sigma ^{2k}}\text{ ,}
\end{equation*}%
and $H_{k}(.)$ denote the well-known Hermite polynomials, which form a
complete orthogonal system in the space $L^{2}\left( \mathbb{R}^{1},(\sqrt{%
2\pi \sigma ^{2}})^{-1}\exp (-z^{2}/2\sigma ^{2}\right) $ of square
integrable functions of Gaussian variables. These polynomials are defined
through the formula: 
\begin{equation*}
H_{j}(z;\sigma ^{2})=(-1)^{j}\sigma ^{2j}\exp \left( \frac{z^{2}}{2\sigma
^{2}}\right) \frac{d^{j}}{dz^{j}}\exp \left( -\frac{z^{2}}{2\sigma ^{2}}%
\right) ,\;\quad j=1,2,\dots
\end{equation*}%
Straightforward computation show that the first five polynomials are: 
\begin{eqnarray}
H_{0}(z) &=&1\text{ , }H_{1}(z)=z\text{ , }H_{2}(z)=z^{2}-\sigma ^{2},\text{ 
}H_{3}(z)=z^{3}-3z\sigma ^{2}\text{ ,}  \notag  \label{hermites} \\
H_{4}(z) &=&z^{4}-6z^{2}\sigma ^{2}+3\sigma ^{4},\text{ }%
H_{5}(z)=z^{5}-10z^{3}\sigma ^{2}+15z\sigma ^{4}\text{ .}  \notag
\end{eqnarray}%
The Hermite polynomials satisfy the differential equation $%
kH_{k-1}(z)=dH_{k}(z)/dz,$ under the boundary conditions $\mathbb{E}%
[H_{m}(Z)]\equiv 0$, where $Z$ is a zero mean Gaussian variable with
variance $\sigma ^{2}$. It is also well-known that, for any mean zero
Gaussian random variables $v$ and $u$, we have: 
\begin{equation}
\mathbb{E}\left[ H_{p}(u)H_{q}(v)\right] =\left\{ 
\begin{array}{rl}
p!\left[ \mathbb{E}(uv)\right] ^{p} & \text{for}\;p=q \\ 
0 & \text{for}\;p\neq q%
\end{array}%
\right. .  \label{eq:ortog}
\end{equation}

The index of the first non null coefficient $b_{k}$ is termed \emph{Hermite
rank }of $g(.)$. Of course, $y_{t}$ is non-Gaussian unless $g(\cdot )$ is
linear.

For our aims, the most important property of Hermite polynomials is their
orthogonality. This property allows us to characterize in a simple way the
dependence structure of a nonlinear transformation of a stationary Gaussian
process that exhibit long range dependence. The following results are
essentially due to \cite{tq75,tq79} and \cite{dm79}. The non-Gaussian case
is more complicated and some results using Appell polynomials are provided
by \cite{sur00}. Let $z_{t}\sim I(d_{z})$; in view of (\ref{eq:covariance})
and (\ref{eq:ortog}) it is easy to see that 
\begin{equation*}
\mathbb{E}[H_{k}(z_{0})H_{k}(z_{\tau })]=k!\gamma ^{k}(\tau )\simeq C\tau
^{k(2d_{z}-1)},\quad \text{as}\;\tau \rightarrow \infty
\end{equation*}%
so the sequence $H_{k+1}(z_{t})$ is ``less''\ dependent then $H_{k}(z_{t})$.
More precisely, if $z_{t}\sim I(d_{z})$, then $H_{k}(z_{t})$ can be viewed
as a long memory series that is fractional integrated of order $d_{k}$ ,%
\begin{equation*}
d_{k}:=\left\{ k\left( d_{z}-\frac{1}{2}\right) +\frac{1}{2}\right\} \vee
0\leq d_{z}\text{ .}
\end{equation*}%
The above equation follows straightforwardly from the equality $%
2d_{k}-1=k(2d_{z}-1)$. The fact that a nonlinear transformation of a
Gaussian process with long range dependence cannot increase its memory can
be used to determine the leading term of the expansion (\ref{eq:hermite}).
For instance, assume for notational simplicity that $z_{t}$ has unit
variance. From (\ref{eq:ortog}) 
\begin{equation*}
\mathbb{E}[H_{k}(z_{t})H_{k}(z_{t+\tau })]=k!\gamma _{zz}^{k}(\tau
)=k!\int_{-\pi }^{\pi }e^{i\lambda \tau }f^{(\ast k)}(\lambda )d\lambda
\end{equation*}%
where 
\begin{equation*}
f^{(\ast k)}(\lambda )=\int_{-\pi }^{\pi }\dots \int_{-\pi }^{\pi
}f_{z}(\lambda -\omega _{1}-\dots -\omega _{k-1})f_{z}(\omega _{1})\dots
f_{z}(\omega _{k-1})d\omega _{1}\dots d\omega _{k-1}
\end{equation*}%
is the $k$ fold convolution of $f_{z}(\lambda )$ with himself (\cite%
{hannan70}, \cite{dalla2004}). The convolution is defined extending $f(.)$
periodically outside of $[-\pi ,\pi ]$. We now have, with an obvious
notation: 
\begin{equation*}
f_{g}(\lambda )=\sum_{k=k_{0}}^{K}b_{k}^{2}k!f_{z}^{(\ast k)}(\lambda )
\end{equation*}%
Therefore, also in the frequency domain, the rank of the transformation
determines the feature of the spectral density, and for $\lambda \rightarrow
0$, the memory of $g(z_{t})$.

Let us now introduce \emph{diagrams}, which are mnemonic devices for
computation of moments and cumulants of polynomial forms in Gaussian random
variables. Our presentation follows \cite{arcones94} and \cite{surgalis03}.
Let $p$ and $\ell _{j}\;,j=1,\dots ,p$, be given integers. A \emph{diagram} $%
G$ of order $(\ell _{1},\dots \ell _{p})$ is a set of points $\{(j,\ell
)):1\leq j\leq p;\;1\leq \ell \leq \ell _{j}\}$ called \emph{vertexes}, and
a set of pairs of these points 
\begin{equation*}
\{((j,\ell ),(k,s)):1\leq j<k<p;1\leq \ell \leq \ell _{j},1\leq s\leq \ell
_{s}\}\text{ ,}
\end{equation*}%
called \emph{edges}. Every vertex is of degree one, that is, it is
considered one time for each graph. We denote by $\mathcal{V}(\ell
_{1},\dots ,\ell _{p})$ the set of diagrams of order $(\ell _{1},\cdots
,\ell _{p})$. The set is empty if $\ell _{1}+\dots +\ell _{p}$ is an odd
number. The set $L_{j}=\{(j,\ell ):1\leq \ell \leq \ell _{j}\}$ is called
the $j$th level of $G$. We will denote the set of edges of diagram $G$ by $%
\mathcal{E}(G)$. Observe that edges connect vertexes of different levels (no
flat edges). A diagram $G$ is said to be \emph{connected} if the rows of the
table cannot be divided in two groups, each of which is partitioned by the
diagram separately. In other words , $G$ is connected if one cannot find a
partition $P_{1}\cup P_{2}=\{1,\dots ,p\}$, $P_{1}\cap P_{2}=\emptyset $, $%
P_{1},P_{2}\neq \emptyset $, such that, for $\mathcal{E}(G)=(g_{1},\dots
g_{k})$, either $g_{i}\in \cup _{j\in P_{1}}L_{j}$ or $g_{i}\in \cup _{j\in
P_{2}}L_{j}$ holds, for $i=1,\dots r$, where $r$ is the number of edges $%
g_{i}$ of the diagram $G$. The set of connected diagrams are indicated by $%
\mathcal{V}^{c}(\ell _{1},\dots ,\ell _{p})$. The main instrument we shall
need below is the following, well known:

\ 

\emph{Diagram Formula:} Let $(z_{1},\dots z_{p})$ be a centered Gaussian
vector, and let $\gamma _{ij}=\mathbb{E}(z_{i}z_{j})$, $i,j=1,\dots ,p$. Let 
$L$ be a table consisting of $p$ rows $\ell _{1},\dots \ell _{p}$, where $%
\ell _{j}$ is the order of Hermite polynomial in the variable $z_{j}$. Then 
\begin{eqnarray*}
\mathbb{E}\left\{ \prod_{j=1}^{p}H_{\ell _{j}}(z_{j})\right\} &=&\sum_{G\in 
\mathcal{V}(\ell _{1},\dots ,\ell _{p})}\prod_{1\leq i<j\leq p}\gamma
_{ij}^{\alpha _{ij}} \\
\text{cum}\left( H_{\ell _{1}}(z_{1}),...,H_{\ell _{p}}(z_{p})\right)
&=&\sum_{G\in \mathcal{V}^{c}(\ell _{1},\dots ,\ell _{p})}\prod_{1\leq
i<j\leq p}\gamma _{ij}^{\alpha _{ij}}
\end{eqnarray*}%
where, for each Gaussian diagram, $\alpha _{ij}$ is the number of edges
between rows $\ell _{i},\ell _{j}$ and cum$\left( H_{\ell
_{1}}(z_{1}),...,H_{\ell _{p}}(z_{p})\right) $ represents the $p-$th order
cumulant.

\ 

Examples of diagrams are represented in Figures 1 to 7 in the Appendix.

\section{Nonlinear cointegration}

We state here more precisely our full set of assumptions.

\ 

\textsc{Assumption A}

1) The following equation holds:%
\begin{equation}
y_{t}=g(x_{t})+e_{t}\text{ },\qquad  \label{eq:recoln}
\end{equation}%
where for $t=1,2,...$ 
\begin{equation*}
\begin{array}{rlll}
g(x_{t})\!\!\!\! & = \sum_{k=k_{0}}^{K}a_{k}x_{t}^{k}\!\!\!\! & =
\sum_{k=k_{0}}^{K}b_{k}H_{k}(x_{t})\text{ , } & b_{k_{0}}\neq 0\text{ },%
\text{ } \\ 
e_{t}\!\!\!\! & = \sum_{\widetilde{k}=\widetilde{k}_{0}}^{\widetilde{K}%
}\theta _{\widetilde{k}}\varepsilon _{t}^{\widetilde{k}} \!\!\!\! & = \sum_{%
\widetilde{k}=\widetilde{k}_{0}}^{\widetilde{K}}\xi _{\widetilde{k}}H_{%
\widetilde{k}}(\varepsilon _{t})\text{ , } & \xi _{\widetilde{k}_{0}}\neq 0%
\text{ },%
\end{array}%
\end{equation*}

2) $(x_{t},\varepsilon _{t})^{\prime }$ are jointly Gaussian and long
memory, that is, as $\tau \rightarrow \infty $%
\begin{eqnarray*}
\gamma _{xx}(\tau ) &\simeq &G_{xx}\tau ^{2d_{x}-1},\text{ }0<G_{xx}<\infty
\\
\gamma _{\varepsilon \varepsilon }(\tau ) &\simeq &G_{\varepsilon
\varepsilon }\tau ^{2d_{\varepsilon }-1},\text{ }0<G_{\varepsilon
\varepsilon }<\infty \\
\gamma _{x\varepsilon }(\tau ) &\simeq &G_{x\varepsilon }\tau
^{d_{x}+d_{\varepsilon }-1},\text{ }|G_{x\varepsilon }|<\infty
\end{eqnarray*}%
for $0\leq d_{\varepsilon },d_{x}<\frac{1}{2}$ .

3) The parameters $K,\widetilde{k}_{0}$ are such that 
\begin{equation*}
K(2d_{x}-1)>\left\{ -1\vee \widetilde{k}_{0}(2d_{\varepsilon }-1)\right\} 
\text{ . }
\end{equation*}

\ 

Assumptions A1-A2 identify a polynomial cointegration model where the
residual is a Gaussian subordinated process. Assumption A3 ensures that $%
H_{K}(x_{t})$ is still a long memory process, with stronger memory than $%
e_{t}.$ This is needed for consistency and indeed it is also a necessary
identification condition: recall $x_{t}$ and $e_{t}$ can be correlated, so
there are no means to distinguish $H_{k}(x_{t})$ and $e_{t}$ unless the
former has stronger long range dependence. Recall that $H_{k}(x_{t})\sim
I(d_{k})$ , $k=k_{0},...,K$, where $2d_{k}-1:=k(2d_{x}-1).$ In this paper,
we take $k_{0}$ and $K$ to be known, whereas their estimation will be
addressed in a different work. Note that to implement our estimates we need
no a priori information on $\widetilde{k}_{0},\widetilde{K}$, although the
value of $\widetilde{k}_{0}(2d_{\varepsilon }-1)$ does affect the rate of
consistency of our estimators.

As mentioned before, (\ref{eq:model}) is a cointegrating relation, so we
allow $\mathbb{E}(x_{t}\varepsilon _{t})$ (and hence $\mathbb{E}(x_{t}e_{t})$%
) to be different from zero. As for linear cointegration, this leads to the
inconsistency of OLS and justifies the use of the spectral regression
methods for the estimation of Hermite coefficients. Concerning the kernel,
we write $k_{M}(\cdot )=k(\tau /M)$ and introduce the following

\ 

\textsc{Assumption }B: The kernel $k(\cdot )$ is a real-valued, symmetric
Lebesgue measurable function that, for $\upsilon \in \mathbb{R}$, satisfies 
\begin{equation*}
\int_{-1}^{1}k(\upsilon )d\upsilon =1\qquad 0\leq k(\upsilon )\leq \infty
,\quad k(\upsilon )=0\quad \text{for}\;|\upsilon |>1.
\end{equation*}

\ 

Our final assumption is a standard bandwidth condition.

\ 

\bigskip \textsc{Assumption }C: Let $\eta =K\vee \widetilde{k}_{0};$ as $%
n\rightarrow \infty $ , 
\begin{equation*}
\frac{1}{M}+\frac{M^{3\vee (\eta -2)}}{n}\rightarrow 0\text{ .}
\end{equation*}

\ 

Assumption C imposes a minimal lower bound and a significant upper bound on
the behaviour of the user-chosen bandwidth parameter $M.$ The need for this
bandwidth condition is made clear by inspection of the proof in the
appendix; heuristically, as $K$ grows the signal in $H_{K}(x_{t})$
decreases, which makes the estimation harder; on the other hand an increase
in $\widetilde{k}_{0}$ makes the convergence rates in Lemma 1 and Theorem 1
faster, whence the need for tighter bandwidth conditions. We are not
claiming Assumption C is sharp, however an inspection of the Proof of Lemma
1 reveals that any improvement is likely to require at least almost
unmanageable computations.

Equation (\ref{eq:recoln}) can be rewritten more compactly as 
\begin{equation*}
y_{t}=\mathbf{\beta }^{\prime }H(x_{t})+e_{t}\text{ , where }%
H(x_{t})=\left\{ H_{1}(x_{t}),...,H_{K}(x_{t})\right\} ^{\prime }\text{ }.
\end{equation*}

Let us now define: 
\begin{equation*}
f_{HH}(\lambda )=\left[ 
\begin{array}{cccc}
f_{11}(\lambda ) & 0 & \cdots  & \cdots  \\ 
0 & f_{22}(\lambda ) & 0 & \cdots  \\ 
\vdots  & \vdots  & \ddots  & \vdots  \\ 
0 & \vdots  & \vdots  & f_{KK}(\lambda )%
\end{array}%
\right] \text{ , }f_{He}(\lambda )=\left[ 
\begin{array}{c}
f_{1e}(\lambda ) \\ 
f_{2e}(\lambda ) \\ 
\vdots  \\ 
f_{Ke}(\lambda )%
\end{array}%
\right] 
\end{equation*}%
and let also, for $a,b=1,2,\dots K$.%
\begin{eqnarray*}
\gamma _{ab}(\tau ) &=&\mathbb{E}\left[ H_{a}(x_{t})H_{b}(x_{t+\tau })\right]
=a!\delta _{a}^{b}\left\{ \mathbb{E}\left( x_{t}x_{t+\tau }\right) \right\}
^{a}, \\
\gamma _{ae}(\tau ) &=&\mathbb{E}\left[ H_{a}(x_{t})e_{t+\tau }\right] =%
\mathbb{E}\left[ H_{a}(x_{t})\sum_{\widetilde{k}=\widetilde{k}_{0}}^{%
\widetilde{K}}\xi _{\widetilde{k}}H_{\widetilde{k}}(\varepsilon _{t})\right] 
\\
&=&\left\{ 
\begin{array}{c}
a!\xi _{a}\left\{ \mathbb{E}\left( x_{t}\varepsilon _{t+\tau }\right)
\right\} ^{a}\text{ for }a\leq \widetilde{K} \\ 
0\text{ , otherwise}%
\end{array}%
\right. .
\end{eqnarray*}%
where $\delta _{a}^{b}$ represents the Kronecker delta function. Likewise%
\begin{eqnarray*}
f_{aa}(\lambda ) &:&=(2\pi )^{-1}\int_{-\infty }^{\infty }\gamma _{aa}(\tau
)e^{-i\lambda \tau }=a!f_{x}^{(\ast a)}(\lambda )\text{ , } \\
f_{ay}(\lambda ) &:&=(2\pi )^{-1}\int_{-\infty }^{\infty }\gamma _{ay}(\tau
)e^{-i\lambda \tau }\text{ , }\gamma _{ay}(\tau ):=\mathbb{E}\left[
H_{a}(x_{t})y_{t+\tau }\right]  \\
f_{ae}(\lambda ) &:&=(2\pi )^{-1}\int_{-\infty }^{\infty }\gamma _{ae}(\tau
)e^{-i\lambda \tau }\text{ .}
\end{eqnarray*}%
The Weighted Covariance Estimator (WCE) of $\beta ^{\prime }=(\beta
_{1},\dots ,\beta _{K})$ is defined as%
\begin{equation*}
\hat{\beta}_{M}=\hat{f}_{HH}(0)^{-1}\hat{f}_{Hy}(0)\text{ ,}
\end{equation*}%
whence%
\begin{equation*}
\hat{\beta}_{M}-\beta =\hat{f}_{HH}(0)^{-1}\hat{f}_{He}(0)\text{ ;}
\end{equation*}%
as usual, we assume $\hat{f}_{HH}(0)$ is non-singular, where 
\begin{equation*}
\hat{f}_{HH}(0)=\frac{1}{2\pi }\left[ 
\begin{array}{ccc}
\sum_{\tau =-M}^{M}k(\tau /M)c_{11}(\tau ) & \cdots  & \sum_{\tau
=-M}^{M}k(\tau /M)c_{1K}(\tau ) \\ 
\vdots  & \ddots  & \vdots  \\ 
\sum_{\tau =-M}^{M}k(\tau /M)c_{K1}(\tau ) & \cdots  & \sum_{\tau
=-M}^{M}k(\tau /M)c_{KK}(\tau )%
\end{array}%
\right] \text{ ,}
\end{equation*}%
\begin{equation*}
\hat{f}_{Hz}(0)=\frac{1}{2\pi }\left[ 
\begin{array}{c}
\sum_{\tau =-M}^{M}k(\tau /M)c_{1z}(\tau ) \\ 
\vdots  \\ 
\sum_{\tau =-M}^{M}k(M)c_{Kz}(\tau )%
\end{array}%
\right] 
\end{equation*}%
and 
\begin{equation*}
\begin{array}{lll}
c_{ab}(\tau )= & \left\{ 
\begin{array}{l}
n^{-1}\sum_{t=1}^{n-\tau }H_{a}(x_{t})H_{b}(x_{t+\tau }) \\ 
n^{-1}\sum_{t=|\tau |+1}^{n}H_{a}(x_{t})H_{b}(x_{t-|\tau |})%
\end{array}%
\right.  & 
\begin{array}{l}
\qquad \tau \geq 0 \\ 
\qquad \tau <0%
\end{array}
\\ 
c_{az}(\tau )= & \left\{ 
\begin{array}{l}
n^{-1}\sum_{t=1}^{n-\tau }H_{a}(x_{t})z_{t+\tau } \\ 
n^{-1}\sum_{t=|\tau |+1}^{n}H_{a}(x_{t})z_{t-|\tau |}%
\end{array}%
\right.  & 
\begin{array}{l}
\qquad \tau \geq 0 \\ 
\qquad \tau <0%
\end{array}%
\end{array}%
\end{equation*}%
for $a,b=1,2,\dots K$, $z=e,y.$ The following lemma is the main tool for our
consistency result, compare Lemma 1 in \cite{mar00}. As before, we write 
\begin{equation*}
d_{a}:=a\left( d_{x}-\frac{1}{2}\right) +\frac{1}{2}\text{ , }d_{e}=\left\{ 
\widetilde{k}_{0}\left( d_{\varepsilon }-\frac{1}{2}\right) +\frac{1}{2}%
\right\} \vee 0\text{ ;}
\end{equation*}%
by Assumption A3 we have $d_{a}>0,$ $a=k_{0},...,K.$\newline

\ \newline
\textbf{LEMMA 1} Under Assumptions A-C, as $n\rightarrow \infty $ we have: 
\begin{align}
\sum_{\tau =-M}^{M}k\left( \frac{\tau }{M}\right) \left\{ c_{ab}(\tau
)-\gamma _{ab}(\tau )\right\} & =o_{p}( M^{d_{a}+d_{b}})  \label{lemma11} \\
\sum_{\tau =-M}^{M}k\left( \frac{\tau }{M}\right) \left\{ c_{ae}(\tau
)-\gamma _{ae}(\tau )\right\} & =o_{p}( M^{d_{a}+d_{e}})  \label{lemma12}
\end{align}
for $a,b=1,2,\dots K$\newline
\textbf{Proof \ }See Appendix

\ 

We are now ready to state the main result of this paper. Let%
\begin{eqnarray*}
B_{ab} &:&=a!G_{xx}^{a}\delta _{a}^{b}\int_{-1}^{1}k(\upsilon )|\upsilon
|^{a(2d_{x}-1)}d\upsilon <\infty \text{ ,} \\
B_{ae} &:&=a!\xi _{a}\left\{ G_{x\varepsilon }\right\}
^{a}\int_{-1}^{1}k(\upsilon )|\upsilon |^{a(d_{x}+d_{\varepsilon
}-1)}d\upsilon <\infty \text{ , for }a\leq \widetilde{K}\text{ ,}
\end{eqnarray*}%
see also Assumption B, $a,b=k_{0},...,K$. Let%
\begin{equation*}
\mathcal{B}_{HH}=\text{diag}\left\{ B_{11},\dots B_{KK}\right\} \text{ , }%
\mathcal{B}_{He}=\left\{ B_{1e},...,B_{Ke}\right\} \text{ , }\mathcal{M}=%
\text{diag}\left\{ M^{-d_{1}},\dots M^{-d_{K}}\right\} \text{ .}
\end{equation*}%
Note that $B_{ae}=0$ unless $a\leq \widetilde{K},$ due to the orthogonality
of Hermite polynomials.

\ \newline
\textbf{Theorem 1 }Under the Assumptions A-C, as $n\rightarrow \infty $ 
\begin{equation*}
\left[ 
\begin{array}{ccc}
M^{d_{1}-d_{e}} & 0 & 0 \\ 
0 & \ddots & 0 \\ 
0 & 0 & M^{d_{K}-d_{e}}%
\end{array}%
\right] \left( \hat{\beta}_{M}-\beta \right) =\mathcal{B}_{HH}^{-1}\mathcal{B%
}_{He}+o_{p}(1)\text{ .}
\end{equation*}%
\textbf{Proof }By the dominated convergence theorem, as $M\rightarrow \infty 
$ 
\begin{align*}
M^{-(d_{a}+d_{b})}\sum_{\tau =-M}^{M}k\left( \frac{\tau }{M}\right) \gamma
_{ab}(\tau )& =\sum_{\tau =-M}^{M}k\left( \frac{\tau }{M}\right) \frac{%
\gamma _{ab}(\tau )}{M^{d_{a}+d_{b}-1}}\frac{1}{M}\rightarrow B_{ab} \\
M^{-(d_{a}+d_{e})}\sum_{\tau =-M}^{M}k\left( \frac{\tau }{M}\right) \gamma
_{ae}(\tau )& =\sum_{\tau =-M}^{M}k\left( \frac{\tau }{M}\right) \frac{%
\gamma _{1e}(\tau )}{M^{d_{a}+d_{e}-1}}\frac{1}{M}\rightarrow B_{ae}
\end{align*}

From Lemma 1, it follows easily that 
\begin{equation*}
\hat{f}_{HH}(0)=\left[ 
\begin{array}{cccc}
\zeta _{1}+o_{p}(M^{2d_{1}}) & o_{p}(M^{d_{1}+d_{2}}) & \cdots & 
o_{p}(M^{d_{1}+d_{p}}) \\ 
o_{p}(M^{d_{2}+d_{1}}) & \zeta _{2}+o_{p}(M^{2d_{2}}) & \cdots & 
o_{p}(M^{d_{2}+d_{p}}) \\ 
\vdots & \vdots & \ddots & \vdots \\ 
o_{p}(M^{d_{K}+d_{1}}) & \cdots & \cdots & \zeta _{K}+o_{p}(M^{2d_{K}})%
\end{array}%
\right]
\end{equation*}%
where%
\begin{equation*}
\zeta _{a}:=\frac{1}{2\pi }\sum_{\tau =-M}^{M}k(\frac{\tau }{M})\gamma
_{aa}(\tau )\text{ .}
\end{equation*}%
Moreover 
\begin{equation*}
\mathcal{M}\hat{f}_{HH}(0)\mathcal{M}=\left[ 
\begin{array}{ccc}
B_{11}+o_{p}(1) & \cdots & o_{p}(1) \\ 
\vdots & \ddots & \vdots \\ 
o_{p}(1) & \cdots & B_{KK}+o_{p}(1)%
\end{array}%
\right] \rightarrow \mathcal{B}_{HH}\text{ .}
\end{equation*}%
Therefore, for $M\rightarrow \infty $ 
\begin{equation*}
\hat{f}_{HH}(0)=\mathcal{M}^{-1}\mathcal{B}_{HH}\mathcal{M}^{-1}+o_{p}(1)=%
\mathcal{B}_{HH}\mathcal{M}^{-2}+o_{p}(1)\text{ ,}
\end{equation*}%
since $\mathcal{B}_{HH}$ is diagonal and hence commutes with $\mathcal{M}%
^{-1}.$ Using the same arguments, it follows easily that: 
\begin{equation*}
M^{-d_{e}}\mathcal{M}\hat{f}_{he}(0)\rightarrow \mathcal{B}_{He}\text{ , as }%
n\rightarrow \infty \text{ .}
\end{equation*}%
Finally, as$\;n\rightarrow \infty $ , 
\begin{equation*}
M^{-d_{e}}\mathcal{M}^{-1}\left\{ \hat{\beta}_{M}-\beta \right\} =\left\{ 
\mathcal{M}\hat{f}_{hh}(0)\mathcal{M}\right\} ^{-1}\mathcal{M}M^{-d_{e}}\hat{%
f}_{he}(0)\rightarrow \mathcal{B}_{HH}^{-1}\mathcal{B}_{He}\text{ ,}\qquad
\end{equation*}%
which completes the proof of Theorem 1.

\hfill$\square$

\ \newline
\textbf{Remark} In Theorem 1 we have proved the consistency of the WCE\
estimator of the cointegrating vector, $\hat{\beta}_{M}\overset{p}{%
\rightarrow }\beta $. In a very loose sense, this result follows from
consistency of a continuously averaged estimate of the spectral density at
frequency zero, see Lemma 1. It is also possible to use Lemma 1 to derive a
robust estimate for the memory parameter of an observed, Gaussian
subordinated series $w_{t}:=g(x_{t}),$ ($k_{0}(d_{x}-\frac{1}{2})+\frac{1}{2}%
=:d_{w},$ say)$.$ We use a very similar idea to the averaged periodogram
estimate advocated by \cite{rob94}. More precisely, with an obvious notation
we can consider%
\begin{eqnarray*}
\widetilde{d}_{w} &:&=\frac{\log \left| \sum_{\tau =-M}^{M}k(\frac{\tau }{M}%
)c_{ww}(\tau )\right| }{2\log M}=d_{w}+\frac{\log B_{ww}}{2\log M}+o_{p}(1)%
\text{ ,} \\
&=&d_{w}+o_{p}(1)\text{ ,}
\end{eqnarray*}%
where we have used Lemma 1. This estimate converges at a mere logarithmic
rate and it is not asymptotically centered around zero; it is however
consistent under much broader circumstances than usually allowed for in the
literature. See also \cite{dalla2004} for very general results on
consistency for long memory estimates.

\section{Comments and conclusions}

We view this paper as a first step in a new research direction, and as such
we are well aware that it leaves several questions unresolved and open for
future research. A first issue relates to the choice of the Hermite rank $%
k_{0}$ and of $K$. As far as the former is concerned, we remark that for the
great majority of practical applications, $k_{0}$ can be taken a priori as 1
or 2. Under the assumption that $k_{0}=1$, the equality $d_{x}=d_{y}$ holds;
this trivial observation immediately suggests a naive test for $k_{0}=1$,
which can be simply implemented by testing for equality of the two memory
parameters. It should be noted, however, that when $x_{t}$ and $y_{t}$ are
cointegrated the standard asymptotic results on multivariate long memory
estimation (for instance \cite{rob95lp}) do not hold. Incidentally, we note
that the nonlinear framework allows to cover the possibility of
cointegration among time series with different integration orders, a
significant extension over the standard paradigm.

For $K$, we can take as an identifying assumption 
\begin{equation}
K:=\text{argmax}(k:k(2d_{x}-1)>(2d_{e}-1))\text{ ;}  \label{idecon}
\end{equation}%
higher order terms can be thought of as included by definition in the
residuals, to make identification possible. Indeed, it is natural to suggest
to view $g(.)$ as a general nonlinear function and envisage $K$ as growing
with $n;$ we expect, however, that only the projection coefficients $b_{k}$
with $k$ satisfying (\ref{idecon}) could be consistently estimated in this
broader framework. On the other hand, we note that the it is also possible
to estimate consistently $K^{\ast }<K$ regression coefficients, by simply
dropping the higher order regressors: it is immediate to see that their
inclusion in the residual would not alter any of our asymptotic result
(there may be an effect in finite samples, however). We stress that a lower
number of regressors allows in general a weaker bandwidth condition, see
Assumption C.

The extension to multivariate regressors does not seem to pose any new
theoretical problem: multivariate generalizations of Hermite expansions are
well known to the literature. Of course, much more challenging seems to be
the possibility to allow for multiple cointegrating relationships. An
important point to remark is the following. In standard cointegration
theory, the role of the variables on the left and on on the right-hand sides
is, by all means, symmetric: this is no longer the case when nonlinear
relationships are allowed. In particular, it should be noted that the memory
parameter of the dependent variable $y_{t}$ is always smaller or equal than $%
d_{x}$; this information can be exploited in an obvious way to decide the
form of the regression, provided that first step estimates of the long
memory parameters are available. We also remark that our procedure requires
a preliminary knowledge on the variance of the regressor $x_{t}$; such
knowledge can clearly be derived from first step estimates, and we leave for
future research the analysis of its consequences in finite samples.

In this paper, we restricted ourselves to consistency results, and gave no
hint on asymptotic distributions. The latter are likely to be non-Gaussian,
at least if the Hermite rank is larger than one and/or the memory of the raw
series is such to make their autocovariances not square summable (see for
instance \cite{fox_taqqu85,fox_taqqu86}. A much wider issue relates to the
possible extension to nonstationary circumstances. Here, a major technical
difficulty arises: the higher order terms in Hermite expansions need no
longer be of smaller order in the presence of nonstationarity. We believe,
however, that the stationary framework considered in this paper is of
sufficient interest by itself for applications to real data, see again \cite%
{cn04} for examples on how fractional cointegration among stationary
variables may be implied by some models of volatility, based on the
Black-Scholes formula for option pricing.

\section*{Appendix}

\textbf{Proof of Lemma 1 }Recall we have 
\begin{equation*}
\gamma _{ab}(\tau )\simeq G|\tau |^{d_{a}+d_{b}-1}\quad \text{ as }\quad
\tau \rightarrow \infty \text{ ,}
\end{equation*}%
where for $a,b=1,...,K,$ $d_{a}$ is such that 
\begin{equation}
d_{a}:=\left\{ 
\begin{array}{cl}
\frac{a}{2}(2d_{x}-1)+\frac{1}{2} & \;\qquad \text{for}\quad \;a(2d_{x}-1)>-1
\\ 
0 & \;\qquad \text{for}\quad \;a(2d_{x}-1)<-1%
\end{array}%
\right. .  \label{eq:index}
\end{equation}%
The first part of the proof follows closely \cite{mar00}. For (\ref{lemma11}%
), it is sufficient to show that 
\begin{eqnarray*}
Var\left\{ \sum_{\tau =-M+1}^{M-1}k\left( \frac{\tau }{M}\right) c_{ab}(\tau
)\right\} &=&\mathbb{E}\left\{ \sum_{\tau =-M+1}^{M-1}k\left( \frac{\tau }{M}%
\right) \left[ c_{ab}(\tau )-\left( 1-\frac{\tau }{n}\right) \gamma
_{ab}(\tau )\right] \right\} ^{2} \\
&\leq
&C\sum_{p=-M}^{M}\sum_{q=-M}^{M}|Cov\{c_{ab}(p),c_{ab}(q)%
\}|=o(M^{2d_{a}+2d_{b}})
\end{eqnarray*}%
From \cite{hannan70}, p.210 we have: 
\begin{eqnarray}
&&Cov\{c_{ab}(p),c_{ab}(q)\}  \notag \\
&=&\!\!\frac{1}{n}\!\sum_{r=-n+1}^{n-1}\!\!\left( 1-\frac{|r|}{n}\!\right)
\left\{ \gamma _{aa}(r)\gamma _{bb}(r\!+\!q\!-\!p)+\gamma
_{ab}(r\!+\!q)\gamma _{ba}(r\!-\!p)\right\}  \label{eq:ones} \\
&&+\frac{1}{n^{2}}\sum_{r=-n+1}^{n-1}\sum_{s=1-r}^{n-r}\text{cum}%
_{abab}\left( s,s+p,s+r,s+r+q\right) \text{ ,}  \label{eq:two}
\end{eqnarray}%
where%
\begin{equation*}
\text{cum}_{abab}\left( s,s+p,s+r,s+r+q\right) =\text{cum}\left\{
H_{a}(x_{s}),H_{b}(x_{s+p}),H_{a}(x_{s+r}),H_{b}(x_{s+r+q})\right\} \text{ .}
\end{equation*}%
Likewise, for (\ref{lemma12}) we shall show that 
\begin{eqnarray*}
Var\left\{ \sum_{\tau =-M+1}^{M-1}k\left( \frac{\tau }{M}\right)
c_{ae}(p)\right\} \!\!\!\! &=&\!\!\!\!\mathbb{E}\left\{ \sum_{\tau
=-M+1}^{M-1}k\left( \frac{\tau }{M}\right) \left[ c_{ae}(\tau )-\left( 1-%
\frac{\tau }{n}\right) \gamma _{ae}(\tau )\right] \right\} ^{2} \\
\!\!\!\! &\leq &\!\!\!\!C\sum_{p=-M}^{M}\sum_{q=-M}^{M}\left|
Cov\{c_{ae}(p),c_{ae}(q)\}\right| =o(M^{2d_{a}+2d_{e}})
\end{eqnarray*}%
where 
\begin{eqnarray}
&&Cov\left\{ c_{ae}(p),c_{ae}(q)\right\}  \notag \\
&=&\!\!\frac{1}{n}\!\sum_{r=-n+1}^{n-1}\!\!\left( 1-\frac{|r|}{n}\right)
\left\{ \gamma _{aa}(r)\gamma _{ee}(r\!+\!q\!-\!p)+\gamma
_{ae}(r\!+\!q)\gamma _{ea}(r\!-\!p)\right\}  \label{mix1} \\
&&+\frac{1}{n^{2}}\sum_{r=-n+1}^{n-1}\sum_{s=1-r}^{n-r}\text{cum}%
_{aeae}\left( s,s+p,s+r,s+r+q\right) \text{ ,}  \label{mix2}
\end{eqnarray}%
and%
\begin{eqnarray*}
&&\text{cum}_{aeae}\left( s,s+p,s+r,s+r+q\right) \\
&=&\text{cum}\left\{ H_{a}(x_{s}),e_{s+p},H_{a}(x_{s+r}),e_{s+r+q}\right\} \\
&=&\sum_{k=\widetilde{k}_{0}}^{\widetilde{K}}\sum_{k^{\prime }=\widetilde{k}%
_{0}}^{\widetilde{K}}\text{cum}\left\{ H_{a}(x_{s}),H_{k}(\varepsilon
_{s+p})H_{a}(x_{s+r}),H_{k^{\prime }},(\varepsilon _{s+r+q})\right\} \text{ .%
}
\end{eqnarray*}%
The argument for (\ref{eq:ones}) and (\ref{mix1}) is the same. For instance,
for (\ref{eq:ones}) we have 
\begin{eqnarray*}
&&\sum_{p=-M}^{M}\sum_{q=-M}^{M}\frac{1}{n}\left| \sum_{r=-n+1}^{n-1}\left(
1-\frac{|r|}{n}\right) \{\gamma _{aa}(r)\gamma _{bb}(r+q-p)\}\right| \\
&\leq &C\frac{M}{n}\sum_{\tau =-2M}^{2M}\left( \sum_{|r|\leq
2M}(|r|+1)^{2d_{a}-1}(|r+\tau |+1)^{2d_{b}-1}+\right. \\
&&
\end{eqnarray*}%
\begin{eqnarray*}
&&+\left. \sum_{|r|>2M}(|r|+1)^{2d_{a}-1}(|r+\tau |+1)^{2d_{b}-1}\right) \\
&=&C\frac{M}{n}\left[ \sum_{|r|\leq 2M}\left( (|r|+1)^{2d_{a}-1}\sum_{\tau
=-2M}^{2M}(|r+\tau |+1)^{2d_{b}-1}\right) \right. \\
&&+\left. \sum_{\tau =-2M}^{2M}\left(
\sum_{2M<|r|<n}(|r|+1)^{2d_{a}-1}(|r+\tau |+1)^{2d_{b}-1}\right) \right] \\
&=&O(Mn^{-1}M^{2d_{a}}M^{2d_{b}})+O(M^{2}n^{-1}n^{2d_{a}+2d_{b}-1})=o(M^{2d_{a}+2d_{b}})%
\text{ .}
\end{eqnarray*}%
As usual, summations over empty sets are taken to be equal to zero. For the
second term we have: 
\begin{eqnarray*}
&&\sum_{p=-M}^{M}\sum_{q=-M}^{M}\frac{1}{n}\left| \sum_{r=-n+1}^{n-1}\left(
1-\frac{|r|}{n}\right) \gamma _{ab}(r+p)\gamma _{ba}(r-q)\right| \\
&\leq &C\sum_{p=-M}^{M}\sum_{q=-M}^{M}\frac{1}{n}\sum_{r=-n+1}^{n-1}\left( 1-%
\frac{|r|}{n}\right) \frac{1}{2}\left| \gamma _{ab}^{2}(r+p)+\gamma
_{ba}^{2}(r-q)\right| \\
&\leq &\frac{C}{n}\sum_{|r|\leq 2M}\left[
\sum_{p=-M}^{M}(|r+p|+1)^{2d_{a}+2d_{b}-2}+%
\sum_{q=-M}^{M}(|r-q|+1)^{2d_{a}+2d_{b}-2}\right] \\
&&+C\frac{M^{2}}{n}\sum_{2M<|r|<n}\left[
(|r+p|+1)^{2d_{a}+2d_{b}-2}+(|r-q|+1)^{2d_{a}+2d_{b}-2}\right] \\
&=&O(Mn^{-1}M^{2d_{a}+2d_{b}})+O(M^{2}n^{-1}n^{2d_{a}+2d_{b}-1})=o(M^{2d_{a}+2d_{b}})%
\text{ .}
\end{eqnarray*}%
The orders of magnitude of the cumulants are investigated by means of the
diagram formula. The proof is quite tedious. We shall focus on (\ref{mix2})
as the argument for (\ref{eq:two}) is entirely analogous. From the diagram
formula it follows easily that, for any finite $k,k^{\prime }\geq \widetilde{%
k}_{0}$ 
\begin{eqnarray}
&&\sum_{k=\widetilde{k}_{0}}^{\widetilde{K}}\sum_{k^{\prime }=\widetilde{k}%
_{0}}^{\widetilde{K}}\left| \text{cum}\left\{ H_{a}(x_{s}),H_{k}(\varepsilon
_{s+p}),H_{a}(x_{s+r}),H_{k^{\prime }}(\varepsilon _{s+r+q})\right\} \right|
\notag \\
&\leq &C\left| \text{cum}\left\{ H_{a}(x_{s}),H_{\widetilde{k}%
_{0}}(\varepsilon _{s+p}),H_{a}(x_{s+r}),H_{\widetilde{k}_{0}}(\varepsilon
_{s+r+q})\right\} \right| .  \label{cumbound}
\end{eqnarray}%
Indeed, increasing the value of $\widetilde{k}_{0}$ to $k,k^{\prime }$
entails including more products of covariances in the cumulant, and these
covariances are bounded. In order to simplify the presentation, we divide it
in three parts, that is

$1)$ $a=1,$ $\widetilde{k}_{0}\geq 2$ or $a\geq 2,$ $\widetilde{k}_{0}=1$

$2)$ $a=2,$ $\widetilde{k}_{0}\geq 2$ or $a\geq 2,$ $\widetilde{k}_{0}=2$

$3)$ $a,\widetilde{k}_{0}\geq 3.$

Throughout the proof, we shall assume for brevity's sake $\widetilde{k}%
_{0}(2d_{\varepsilon }-1)>-1;$ it is simple to check that for $\widetilde{k}%
_{0}(2d_{\varepsilon }-1)\leq -1$ the proof is analogous, indeed slightly
simpler. \ 

\subsubsection*{Part I: $a=1,$ $\widetilde{k}_{0}\geq 2$ or $a\geq 2,$ $%
\widetilde{k}_{0}=1$}

For $a=1,$ $\widetilde{k}_{0}=2$ we have%
\begin{eqnarray*}
&&\sum_{p=-M}^{M}\sum_{q=-M}^{M}\frac{1}{n^{2}}\left| \text{cum}\left\{
x_{s},H_{2}(\varepsilon _{s+p}),x_{s+r},H_{2}(\varepsilon _{s+r+q})\right\}
\right| \\
&\leq &\sum_{p=-M}^{M}\sum_{q=-M}^{M}\frac{C}{n^{2}}\left|
\sum_{r=-n+1}^{n-1}\sum_{s=1-r}^{n-r}\gamma _{x\varepsilon }(p)\gamma
_{x\varepsilon }(q)\gamma _{\varepsilon \varepsilon }(r+q-p)\right. \\
&&\left. +\gamma _{\varepsilon \varepsilon }(r+q-p)\gamma _{\varepsilon
x}(r-p)\gamma _{x\varepsilon }(r+q)\right| \\
&\leq &\frac{C}{n}\sum_{p=-M}^{M}(|p|+1)^{d_{x}+d_{\varepsilon
}-1}\sum_{q=-M}^{M}(|q|+1)^{d_{x}+d_{\varepsilon }-1}\left( \sum_{|r|\leq
3M}(|r+q-p|+1)^{2d_{\varepsilon }-1}\right. \\
&&+\left. \sum_{3M<|r|\leq n}(|r+q-p|+1)^{2d_{\varepsilon }-1}\right) \\
&&+\frac{C}{n}\!\!\sum_{p=-M}^{M}\!\sum_{q=-M}^{M}\!\!\left( \sum_{|r|\leq
3M}(|r\!+\!q\!-\!p|\!+\!1)^{2d_{\varepsilon
}-1}(|r\!-\!p|\!+\!1)^{d_{x}+d_{\varepsilon
}-1}(|r\!+\!q|+1)^{d_{x}+d_{\varepsilon }-1}\right. \\
&&+\left. \sum_{3M<|r|\leq n}(|r+q-p|+1)^{2d_{\varepsilon
}-1}(|r-p|+1)^{d_{x}+d_{\varepsilon }-1}(|r+q|+1)^{d_{x}+d_{\varepsilon
}-1}\right) \\
&=&O(n^{-1}M^{d_{x}+d_{\varepsilon }}M^{d_{x}+d_{\varepsilon
}}M^{2d_{\varepsilon }})+O(n^{-1}M^{2d_{x}+2d_{\varepsilon
}}n^{2d_{\varepsilon }}) \\
&&+O(n^{-1}MM^{d_{x}+d_{\varepsilon }}M^{2d_{\varepsilon
}})+O(n^{-1}M^{2}n^{2d_{x}+4d_{\varepsilon }-3}) \\
&=&O\left( \frac{M}{n}M^{2d_{x}+4d_{\varepsilon }-1}\right)
+o(M^{4d_{\varepsilon }+2d_{x}})+O\left( \frac{M^{2}}{n}M^{d_{x}+3d_{%
\varepsilon }-1})+o(M^{4d_{\varepsilon }+2d_{x}}\right) \\
&=&o(M^{2d_{x}+4d_{\varepsilon }-1})=o(M^{2d_{x}+2d_{e}})\quad \text{
because }2d_{e}=4d_{\varepsilon }-1\text{ .}
\end{eqnarray*}%
The extension to $\widetilde{k}_{0}>2$ is trivial: {%
\begin{eqnarray*}
&&\text{cum}\left\{ x_{s},H_{\widetilde{k}_{0}}(\varepsilon
_{s+p}),x_{s+r},H_{\widetilde{k}_{0}}(\varepsilon _{s+r+q})\right\} \\
&=&\sum_{p=-M}^{M}\sum_{q=-M}^{M}\frac{C}{n^{2}}\Bigg|\sum_{r=-n+1}^{n-1}%
\sum_{s=1-r}^{n-r}\gamma _{x\varepsilon }(p)\gamma _{x\varepsilon }(q)\gamma
_{\varepsilon \varepsilon }^{\widetilde{k}_{0}-1}(r+q-p)
\end{eqnarray*}%
\begin{eqnarray*}
&&+\gamma _{\varepsilon \varepsilon }^{\widetilde{k}_{0}-1}(r+q-p)\gamma
_{\varepsilon x}(r-p)\gamma _{x\varepsilon }(r+q)\Bigg| \\
&=&O(n^{-1}M^{2d_{x}+2d_{\varepsilon }}M^{\widetilde{(k}_{0}-1)(2d_{%
\varepsilon }-1)+1})+O(n^{-1}M^{2d_{x}+2d_{\varepsilon }}n^{\widetilde{(k}%
_{0}-1)(2d_{\varepsilon }-1)+1}) \\
&&+O(n^{-1}M^{d_{x}+d_{\varepsilon }+1}M^{\widetilde{(k}_{0}-1)(2d_{%
\varepsilon }-1)+1})+O(n^{-1}M^{2}n^{2d_{x}+2d_{\varepsilon }-2+\widetilde{(k%
}_{0}-1)(2d_{\varepsilon }-1)+1}) \\
&=&O\left( \frac{M}{n}M^{2d_{x}+2d_{\varepsilon }-1}M^{\widetilde{(k}%
_{0}-1)(2d_{\varepsilon }-1)+1}\right) +o(M^{2d_{x}+\widetilde{k}%
_{0}(2d_{\varepsilon }-1)+1}) \\
&&+O\left( \frac{M^{2}}{n}M^{d_{x}+d_{\varepsilon }-1}M^{\widetilde{(k}%
_{0}-1)(2d_{\varepsilon }-1)+1}\right) +O(n^{-1}M^{2}n^{2d_{x}+\widetilde{k}%
_{0}(2d_{\varepsilon }-1)}) \\
&=&o(M^{2d_{x}+2d_{e}}),
\end{eqnarray*}%
by the same argument as before. The proof for }$a\geq 2,$ $\widetilde{k}%
_{0}=1$ is entirely analogous and hence omitted.

\subsubsection*{Part II: $a=2,$ $\widetilde{k}_{0}\geq 2$ or $a\geq 2,$ $%
\widetilde{k}_{0}=2$}

For $a=2,$ $\widetilde{k}_{0}=2$ we have%
\begin{eqnarray*}
&&\sum_{p=-M}^{M}\sum_{q=-M}^{M}\frac{1}{n^{2}}\left|
\sum_{r=-n+1}^{n-1}\sum_{s=1-r}^{n-r}\text{cum}\{H_{2}(x_{s})H_{2}(%
\varepsilon _{s+p})H_{2}(x_{s+r})H_{2}(\varepsilon _{s+r+q})\}\right| \\
&\leq &\sum_{p=-M}^{M}\sum_{q=-M}^{M}\frac{C}{n^{2}}\Bigg|%
\sum_{r=-n+1}^{n-1}\sum_{s=1-r}^{n-r}\gamma _{x\varepsilon }(p)\gamma
_{x\varepsilon }(q)\gamma _{xx}(r)\gamma _{\varepsilon \varepsilon }(r+q-p)
\\
&&+\gamma _{x\varepsilon }(p)\gamma _{x\varepsilon }(r+q)\gamma
_{\varepsilon x}(r-p)\gamma _{x\varepsilon }(q) \\
&&+\gamma _{xx}(r)\gamma _{\varepsilon \varepsilon }(r+q-p)\gamma
_{\varepsilon x}(r-p)\gamma _{x\varepsilon }(r+q)\Bigg| \\
&\leq &\!\!\!\frac{C}{n}\!\left\{
\sum_{p=-M}^{M}\sum_{q=-M}^{M}\sum_{|r|\leq 3M}\!\left[ (|p|\!+%
\!1)^{d_{x}+d_{\varepsilon }-1}(|q|\!+\!1)^{d_{x}+d_{\varepsilon
}-1}(|r|\!+\!1)^{2d_{x}-1}(|r\!+\!q\!-\!p|\!+\!1)^{2d_{\varepsilon
}-1}\right. \right. \\
&&+(|p|\!+\!1)^{d_{x}+d_{\varepsilon }-1}(|q|\!+\!1)^{d_{x}+d_{\varepsilon
}-1}(|r\!+\!q|\!+\!1)^{d_{x}+d_{\varepsilon
}-1}(|r-p|+1)^{d_{x}+d_{\varepsilon }-1} \\
&&+\left. \left.
(|r|\!+\!1)^{2d_{x}-1}(|r\!+\!q\!-\!p|\!+\!1)^{2d_{\varepsilon
}-1}(|r\!-\!p|\!+\!1)^{d_{x}+d_{\varepsilon
}-1}(|r\!+\!q|\!+\!1)^{d_{x}+d_{\varepsilon }-1}\right] \right\} \\
&&+\frac{C}{n}\!\left\{ \sum_{p=-M}^{M}\sum_{q=-M}^{M}\sum_{3M<r\leq n}\!\!%
\left[ (|p|\!+\!1)^{d_{x}+d_{\varepsilon
}-1}(|q|\!+\!1)^{d_{x}+d_{\varepsilon
}-1}(|r|\!+\!1)^{2d_{x}-1}(|r\!+\!q\!-\!p|+1)^{2d_{\varepsilon }-1}\right.
\right. \\
&&+(|p|\!+\!1)^{d_{x}+d_{\varepsilon }-1}(|q|\!+\!1)^{d_{x}+d_{\varepsilon
}-1}(|r\!+\!q|\!+\!1)^{d_{x}+d_{\varepsilon
}-1}(|r\!-\!p|\!+\!1)^{d_{x}+d_{\varepsilon }-1} \\
&&+\left. \left.
(|r|\!+\!1)^{2d_{x}-1}(|r\!+\!q\!-\!p|\!+\!1)^{2d_{\varepsilon
}-1}(|r\!-\!p|\!+\!1)^{d_{x}+d_{\varepsilon
}-1}(|r\!+\!q|\!+\!1)^{d_{x}+d_{\varepsilon }-1}\right] \right\} \\
&=&O(n^{-1}M^{2d_{x}+2d_{\varepsilon }}M^{2d_{\varepsilon
}})+O(n^{-1}M^{d_{x}+d_{\varepsilon }}M^{d_{x}+d_{\varepsilon
}}M^{d_{x}+d_{\varepsilon }})+O(n^{-1}M^{3d_{\varepsilon }+3d_{x}}) \\
&&+O(n^{-1}M^{2d_{x}+2d_{\varepsilon }}n^{2d_{x}+2d_{\varepsilon
}-1})+O(n^{-1}M^{2d_{x}+2d_{\varepsilon }}n^{2d_{x}+2d_{\varepsilon
}-1})+O(M^{2}n^{-1}n^{4d_{x}+4d\varepsilon -3}) \\
&&
\end{eqnarray*}%
\begin{eqnarray*}
&=&O\left( \frac{M^{2}}{n}M^{2d_{x}+4d_{\varepsilon }-2}\right) +O\left( 
\frac{M^{2}}{n}M^{2d_{x}+4d_{\varepsilon }-2}\right)
+o(M^{4d_{x}+4d_{\varepsilon }-2}) \\
&=&o(M^{2d_{2}+2d_{e}})\text{ because }2d_{2}=4d_{x}-1\text{ and }%
2d_{e}=4d_{\varepsilon }-1\text{ .}
\end{eqnarray*}%
For $\widetilde{k}_{0}>2$ the argument is very much the same:%
\begin{eqnarray*}
&&\sum_{p=-M}^{M}\sum_{q=-M}^{M}\frac{1}{n^{2}}\left|
\sum_{r=-n+1}^{n-1}\sum_{s=1-r}^{n-r}\text{cum}\{H_{2}(x_{s})H_{\widetilde{k}%
_{0}}(\varepsilon _{s+p})H_{2}(x_{s+r})H_{\widetilde{k}_{0}}(\varepsilon
_{s+r+q})\}\right| \\
&\leq &\sum_{p=-M}^{M}\sum_{q=-M}^{M}\frac{C}{n^{2}}\Bigg|%
\sum_{r=-n+1}^{n-1}\sum_{s=1-r}^{n-r}\gamma _{x\varepsilon }(p)\gamma
_{x\varepsilon }(q)\gamma _{xx}(r)\gamma _{\varepsilon \varepsilon }^{%
\widetilde{k}_{0}-1}(r+q-p) \\
&&+\gamma _{x\varepsilon }(p)\gamma _{x\varepsilon }(q)\gamma _{x\varepsilon
}(r+q)\gamma _{\varepsilon x}(r-p)\gamma _{\varepsilon \varepsilon }^{%
\widetilde{k}_{0}-2}(r+q-p) \\
&&+\gamma _{xx}(r)\gamma _{\varepsilon x}(r-p)\gamma _{x\varepsilon
}(r+q)\gamma _{\varepsilon \varepsilon }^{\widetilde{k}_{0}-1}(r+q-p)\Bigg|
\\
&=&O(n^{-1}M^{2d_{x}+2d_{\varepsilon }}M^{\widetilde{(k}_{0}-1)(2d_{%
\varepsilon }-1)+1})+O(n^{-1}M^{2d_{x}+2d_{\varepsilon }}M^{\widetilde{(k}%
_{0}-2)(2d_{\varepsilon }-1)+1}) \\
&&+O(n^{-1}M^{2d_{x}}M^{d_{x}+d_{\varepsilon }}M^{\widetilde{(k}%
_{0}-1)(2d_{\varepsilon }-1)+1})+O(n^{-1}M^{2d_{x}+2d_{\varepsilon
}}n^{2d_{x}-1+\widetilde{(k}_{0}-1)(2d_{\varepsilon }-1)+1}) \\
&&+O(n^{-1}M^{2d_{x}+2d_{\varepsilon }}n^{2d_{x}+2d_{\varepsilon }-1+%
\widetilde{(k}_{0}-2)(2d_{\varepsilon
}-1)})+O(M^{2}n^{-1}n^{4d_{x}+2d_{\varepsilon }-2+\widetilde{(k}%
_{0}-1)(2d_{\varepsilon }-1)}) \\
&=&O\left( \frac{M^{2}}{n}M^{2d_{x}-1}M^{\widetilde{k}_{0}(2d_{\varepsilon
}-1)+1}\right) +o\left( \frac{M^{2}}{n}M^{2d_{x}+2d_{\varepsilon }-1+%
\widetilde{k}_{0}(2d_{\varepsilon }-1)+1}\right) \\
&&+O\left( \frac{M^{2}}{n}M^{3d_{x}-1-d_{\varepsilon }}M^{\widetilde{k}%
_{0}(2d_{\varepsilon }-1)+1}\right) +O\left( \frac{M^{2}}{n}M^{4d_{x}-1}M^{%
\widetilde{k}_{0}(2d_{\varepsilon }-1)+1}\right) \\
&=&o(M^{2d_{2}+2d_{e}})\text{ , because }2d_{e}=\widetilde{k}%
_{0}(2d_{\varepsilon }-1)+1\text{ .}
\end{eqnarray*}

\subsection*{Part III: $a\geq 3,$ $\widetilde{k}_{0}\geq 3$}

We note that, by the diagram formula (as in (\ref{cumbound}))%
\begin{eqnarray*}
&&\left| \text{cum}\left[ H_{a}(x_{s})H_{\widetilde{k}_{0}}(\varepsilon
_{s+p})H_{a}(x_{s+r})H_{\widetilde{k}_{0}}(\varepsilon _{s+r+q})\right]
\right| \\
&\leq &C\left| \text{cum}\left[ H_{3}(x_{s})H_{3}(\varepsilon
_{s+p})H_{3}(x_{s+r})H_{3}(\varepsilon _{s+r+q})\right] \right| \text{ }.
\end{eqnarray*}%
It suffices then to focus on $a=\widetilde{k}_{0}=3.$ There are seven
different kinds of connected diagrams, which are represented in Figures 1 to
7. We have 
\begin{eqnarray*}
&&\sum_{p=-M}^{M}\sum_{q=-M}^{M}\frac{1}{n^{2}}\left|
\sum_{r=-n+1}^{n-1}\sum_{s=1-r}^{n-r}\text{cum}\{H_{3}(x_{s})H_{3}(%
\varepsilon _{s+p})H_{3}(x_{s+r})H_{3}(\varepsilon _{s+r+q})\}\Bigg|\right.
\\
&=&\sum_{p=-M}^{M}\sum_{q=-M}^{M}\frac{C}{n^{2}}\Bigg|\sum_{r=-n+1}^{n-1}%
\sum_{s=1-r}^{n-r}\gamma _{x\varepsilon }^{2}(p)\gamma _{x\varepsilon
}^{2}(q)\gamma _{\varepsilon x}(r-p)\gamma _{x\varepsilon }(r+q)+
\end{eqnarray*}%
\begin{eqnarray*}
&+&\gamma _{x\varepsilon }(p)\gamma _{xx}(r)\gamma _{x\varepsilon
}(r+q)\gamma _{\varepsilon \varepsilon }(r+q-p)\gamma _{\varepsilon
x}(r-p)\gamma _{x\varepsilon }(q) \\
&+&\gamma _{xx}^{2}(r)\gamma _{\varepsilon \varepsilon }^{2}(r+q-p)\gamma
_{\varepsilon x}(r-p)\gamma _{x\varepsilon }(r+q) \\
&+&\gamma _{x\varepsilon }^{2}(p)\gamma _{x\varepsilon }^{2}(q)\gamma
_{xx}(r)\gamma _{\varepsilon \varepsilon }(r+q-p) \\
&+&\gamma _{xx}^{2}(r)\gamma _{\varepsilon \varepsilon }^{2}(r+q-p)\gamma
_{x\varepsilon }(p)\gamma _{x\varepsilon }(q) \\
&+&\gamma _{x\varepsilon }^{2}(r-p)\gamma _{x\varepsilon }^{2}(r+q)\gamma
_{x\varepsilon }(p)\gamma _{x\varepsilon }(q) \\
&+&\left. \gamma _{x\varepsilon }^{2}(r-p)\gamma _{x\varepsilon
}^{2}(r+q)\gamma _{xx}(r)\gamma _{\varepsilon \varepsilon }(r+q-p)\right|
\end{eqnarray*}%
\begin{eqnarray}
&\leq &\frac{C}{n}\sum_{p=-M}^{M}(|p|+1)^{2(d_{x}+d_{\varepsilon
}-1)}\sum_{q=-M}^{M}(|q|+1)^{2(d_{x}+d_{\varepsilon }-1)}  \notag \\
&\times &\left[ \sum_{|r|\leq 2M}(|r-p|+1)^{d_{x}+d_{\varepsilon
}-1}(|r+q|+1)^{d_{x}+d_{\varepsilon }-1}\right.  \notag \\
&+&\left. \sum_{2M<|r|\leq n}(|r-p|+1)^{d_{x}+d_{\varepsilon
}-1}(|r+q|+1)^{d_{x}+d_{\varepsilon }-1}\right]  \label{one}
\end{eqnarray}%
\begin{eqnarray}
&+&\frac{C}{n}\sum_{p=-M}^{M}(|p|+1)^{d_{x}+d_{\varepsilon
}-1}\sum_{q=-M}^{M}(|q|+1)^{d_{x}+d_{\varepsilon }-1}  \notag \\
&&\times \left[ \sum_{|r|\leq
3M}(|r|+1)^{2d_{x}-1}(|r+q-p|+1)^{2d_{\varepsilon
}-1}(|r-p|+1)^{d_{x}+d_{\varepsilon }-1}(|r+q|+1)^{d_{x}+d_{\varepsilon
}-1}\right.  \notag \\
&+&\left. \!\!\!\sum_{3M<|r|\leq
n}(|r|\!+\!1)^{2d_{x}-1}(|r\!+\!q\!-\!p|\!+\!1)^{2d_{\varepsilon
}-1}(|r\!-\!p|\!+\!1)^{d_{x}+d_{\varepsilon
}-1}(|r\!+\!q|\!+\!1)^{d_{x}+d_{\varepsilon }-1}\right]  \label{two}
\end{eqnarray}
\begin{eqnarray}
+ &&\frac{C}{n}\left[ \left( \sum_{|r|\leq
3M}(|r|+1)^{2(2d_{x}-1)}\sum_{p=-M}^{M}\sum_{q=-M}^{M}(|r+q-p|+1)^{2(2d_{%
\varepsilon }-1)}\right. \right.  \notag \\
&&\times \!\left. (|r\!-\!p|\!+\!1)^{d_{x}+d_{\varepsilon
}-1}(|r\!+\!q|\!+\!1)^{d_{x}+d_{\varepsilon }-1}\right)
+\sum_{p=-M}^{M}\sum_{q=-M}^{M}\left( \sum_{3M<|r|\leq
n}(|r|\!+\!1)^{2(2d_{x}-1)}\right.  \notag \\
&&\times \left. \left. |r+q-p|+1)^{2(2d_{\varepsilon
}-1)}(|r-p|+1)^{d_{x}+d_{\varepsilon }-1}(|r+q|+1)^{d_{x}+d_{\varepsilon
}-1}\right) \right]  \label{three}
\end{eqnarray}%
\begin{eqnarray}
+ &&\frac{C}{n}\sum_{p=-M}^{M}(|p|+1)^{2(d_{x}+d_{\varepsilon
}-1)}\sum_{q=-M}^{M}(|q|+1)^{2(d_{x}+d_{\varepsilon }-1)}\left[ \left(
\sum_{|r|\leq 3M}(|r|+1)^{2d_{x}-1}\right. \right.  \notag \\
&&\times \left. (|r\!+\!q\!-\!p|\!+\!1)^{2d_{x}-1}\right) +\!\left.
\!\sum_{3M<|r|\leq
n}(|r|\!+\!1)^{2d_{x}-1}(|r\!+\!q\!-\!p|\!+\!1)^{2d_{x}-1} \right]
\label{four}
\end{eqnarray}%
\begin{eqnarray}
&+&\frac{C}{n}\sum_{p=-M}^{M}(|p|+1)^{d_{x}+d_{\varepsilon
}-1}\sum_{q=-M}^{M}(|q|+1)^{d_{x}+d_{\varepsilon }-1}\left[ \left(
\sum_{|r|\leq 3M}(|r|+1)^{2(2d_{x}-1)}\phantom{cccccccc}\right. \right. 
\notag \\
&&\times \left. (|r\!+\!q\!-\!p|\!+\!1)^{2(2d_{\varepsilon }-1)}\right)
+\left. \sum_{3M<|r|\leq
n}(|r|\!+\!1)^{2(2d_{x}-1)}(|r\!+\!q\!-\!p|+1)^{2(2d_{\varepsilon }-1)} 
\right]  \label{five}
\end{eqnarray}%
\begin{eqnarray}
+ &&\frac{C}{n}\sum_{p=-M}^{M}(|p|+1)^{d_{x}+d_{\varepsilon
}-1}\sum_{q=-M}^{M}(|q|+1)^{d_{x}+d_{\varepsilon }-1}\left[ \left(
\sum_{|r|\leq 3M}(|r+q|+1)^{2(d_{x}+d_{\varepsilon }-1)}\phantom{ccc}\right.
\right.  \notag \\
&&\times \left. \left. (|r\!-\!p|\!+\!1)^{2(d_{x}+d_{\varepsilon
}-1)}\right) +\!\!\!\sum_{3M<|r|\leq
n}(|r\!+\!q|\!+\!1)^{2(d_{x}+d_{\varepsilon
}-1)}(|r\!-\!p|\!+\!1)^{2(d_{x}+d_{\varepsilon }-1)}\right]  \label{six}
\end{eqnarray}%
\begin{eqnarray}
&+&\frac{C}{n}\sum_{|r|\leq
3M}(|r|+1)^{2d_{x}-1}\sum_{p=-M}^{M}\sum_{q=-M}^{M}(|r+p-q|+1)^{2d_{%
\varepsilon }-1}  \notag \\
&&\times (|r-p|+1)^{2(d_{x}+d_{\varepsilon
}-1)}(|r+q|+1)^{2(d_{x}+d_{\varepsilon }-1)}+\sum_{3M<|r|\leq
n}(|r|+1)^{2d_{x}-1}\phantom{ccccccc}  \notag \\
&&\times
\sum_{p=-M}^{M}\sum_{q=-M}^{M}(|r\!+\!p\!-\!q|\!+\!1)^{2d_{\varepsilon
}-1}(|r\!-\!p|\!+\!1)^{2(d_{x}+d_{\varepsilon
}-1)}(|r\!+\!q|\!+\!1)^{2(d_{x}+d_{\varepsilon }-1)}.  \label{seven}
\end{eqnarray}%
After lengthy but straightforward computations, it is not difficult to see
that 
\begin{eqnarray*}
(\ref{one}) &=&O(n^{-1}M^{2d_{x}+2d_{\varepsilon
}-1}M^{2d_{x}+2d_{\varepsilon }-1}M^{d_{x}+d_{\varepsilon
}})+O(n^{-1}M^{4d_{x}-1}M^{4d_{\varepsilon }-1}n^{2d_{x}+2d_{\varepsilon
}-1}) \\
&=&o\left( \frac{M^{4d_{x}+4d_{\varepsilon }-1}}{n}\right)
+O(M^{4d_{x}-1}M^{4d_{\varepsilon }-1}n^{2d_{x}+2d_{\varepsilon }-2}) \\
(\ref{two}) &=&O(n^{-1}M^{d_{x}+d_{\varepsilon }}M^{d_{x}+d_{\varepsilon
}}M^{2d_{x}+2d_{\varepsilon }-1})+O(n^{-1}M^{d_{x}+d_{\varepsilon
}}M^{d_{x}+d_{\varepsilon }}n^{4d_{x}+4d_{\varepsilon }-3}) \\
&=&O\left( \frac{M^{4d_{x}+4d_{\varepsilon }-1}}{n}\right)
+o(M^{4d_{x}-1}M^{4d_{\varepsilon }-1}n^{2d_{x}+2d_{\varepsilon }-2}) \\
(\ref{three}) &=&O(n^{-1}M^{4d_{x}-1}M^{4d_{\varepsilon
}})+O(n^{-1}M^{2}n^{6d_{x}+6d_{\varepsilon }-5}) \\
&=&O\left( \frac{M^{4d_{x}+4d_{\varepsilon }-1}}{n}\right)
+o(M^{4d_{x}-1}M^{4d_{\varepsilon }-1}n^{2d_{x}+2d_{\varepsilon }-2}) \\
(\ref{four}) &=&O(n^{-1}M^{2d_{x}+2d_{\varepsilon
}-1}M^{2d_{x}+2d_{\varepsilon
}-1}M^{2d_{x}})+O(n^{-1}M^{4d_{x}-1}M^{4d_{\varepsilon
}-1}n^{2d_{x}+2d_{\varepsilon }-1}) \\
&=&o\left( \frac{M^{4d_{x}+4d_{\varepsilon }-1}}{n}\right)
+O(M^{4d_{x}-1}M^{4d_{\varepsilon }-1}n^{2d_{x}+2d_{\varepsilon }-2}) \\
(\ref{five}) &=&O(n^{-1}M^{d_{x}+d_{\varepsilon }}M^{d_{x}+d_{\varepsilon
}}M^{4d_{\varepsilon
}-1})+O(n^{-1}M^{2d_{x}}M^{2d_{x}}n^{4d_{x}+4d_{\varepsilon }-3}) \\
&=&O\left( \frac{M^{4d_{x}+4d_{\varepsilon }-1}}{n}\right)
+o(M^{4d_{x}-1}M^{4d_{\varepsilon }-1}n^{2d_{x}+2d_{\varepsilon }-2}) \\
\end{eqnarray*}
\begin{eqnarray*}
(\ref{six}) &=&O(n^{-1}M^{d_{x}+d_{\varepsilon }}M^{d_{x}+d_{\varepsilon
}}M^{2d_{x}+2d_{\varepsilon }-1})+O(n^{-1}M^{2d_{x}+2d_{\varepsilon
}}n^{4d_{x}+4d_{\varepsilon }-3}) \\
&=&O\left( \frac{M^{4d_{x}+4d_{\varepsilon }-1}}{n}\right)
+o(M^{4d_{x}-1}M^{4d_{\varepsilon }-1}n^{2d_{x}+2d_{\varepsilon }-2}) \\
(\ref{seven}) &=&O(n^{-1}M^{2d_{x}}M^{2d_{x}+2d_{\varepsilon
}})+O(n^{-1}M^{2}n^{6d_{x}+6d_{\varepsilon }-5}) \\
&=&o(\frac{M^{4d_{x}+4d_{\varepsilon }-1}}{n})+o(M^{4d_{x}-1}M^{4d_{%
\varepsilon }-1}n^{2d_{x}+2d_{\varepsilon }-2})
\end{eqnarray*}
\begin{equation*}
=O(\frac{M^{4d_{x}+4d_{\varepsilon }-1}}{n})+o(M^{4d_{x}-1}M^{4d_{%
\varepsilon }-1}n^{2d_{x}+2d_{\varepsilon }-2})\text{ .}
\end{equation*}
In view of the previous results, our proof will be completed if we show that 
\begin{equation*}
\frac{M^{4d_{x}+4d_{\varepsilon }-1}}{n}+M^{4d_{x}-1}M^{4d_{\varepsilon
}-1}n^{2d_{x}+2d_{\varepsilon }-2}=o(M^{2d_{a}+2d_{e}})\text{ ,}
\end{equation*}%
where 
\begin{equation*}
2d_{a}+2d_{e}=2ad_{x}+2\widetilde{k}_{0}d_{\varepsilon }-(a+\widetilde{k}%
_{0})+2\text{ .}
\end{equation*}%
We note first that%
\begin{equation*}
\frac{M^{4d_{x}+4d_{\varepsilon }-1}}{nM^{2d_{a}+2d_{e}}}=\frac{%
M^{4d_{x}+4d_{\varepsilon }-1}}{n(M^{2ad_{x}+2\widetilde{k}%
_{0}d_{\varepsilon }-(a+\widetilde{k}_{0})+2})}=\frac{M^{2(2-a)d_{x}+2(2-%
\widetilde{k}_{0})d_{\varepsilon }-3+a+\widetilde{k}_{0}}}{n}\text{ .}
\end{equation*}%
From (\ref{eq:index}) and $d_{e}>0$ it follows that 
\begin{equation*}
d_{x}>\frac{1}{2}-\frac{1}{2a}\quad \text{and}\quad d_{\varepsilon }>\frac{1%
}{2}-\frac{1}{2\widetilde{k}_{0}}\text{ ,}
\end{equation*}%
whence, because $a,\widetilde{k}_{0}\geq 2$ we have 
\begin{eqnarray*}
\frac{M^{2(2-a)d_{x}+2(2-\widetilde{k}_{0})d_{\varepsilon }-3+a+\widetilde{k}%
_{0}}}{n} &\leq &\frac{M^{2(2-a)(\frac{1}{2}-\frac{1}{2a})+2(2-\widetilde{k}%
_{0})(\frac{1}{2}-\frac{1}{2\widetilde{k}_{0}})-3+a+\widetilde{k}_{0}}}{n} \\
&=&o\left( \frac{M^{3}}{n}\right) =o(1)\text{ ,}
\end{eqnarray*}%
in view of Assumption C. To complete the proof, note that, again from
Assumption C, for some $\alpha >a-2,\widetilde{k}_{0}-2$ we have $M^{\alpha
}=O(n),$ whence%
\begin{equation*}
\frac{M^{4d_{x}+4d_{\varepsilon }-2}n^{2d_{x}+2d_{\varepsilon }-2}}{%
M^{2ad_{x}+2\widetilde{k}_{0}d_{\varepsilon }-(a+\widetilde{k}_{0})+2}}%
=o(M^{(2-a+\alpha )(2d_{x}-1)+(2-\widetilde{k}_{0}+\alpha )(2d_{\varepsilon
}-1)})=o(1)\qquad \text{as}\quad n\rightarrow \infty \text{ .}
\end{equation*}%
Thus the proof is completed. \hfill $\square$ \newpage 
\begin{figure}[!h]
\centering \includegraphics{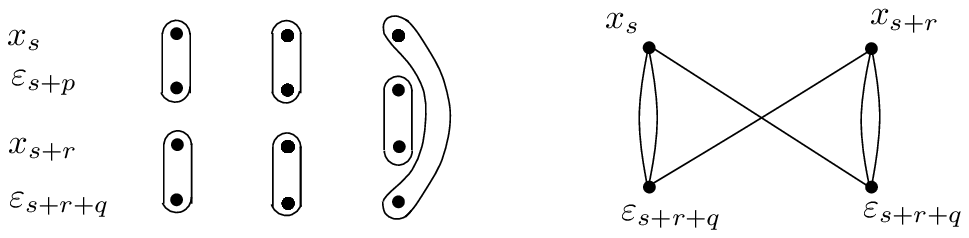} \label{fig:grafo1}
\caption{$\protect\gamma _{x\protect\varepsilon }^{2}(p)\protect\gamma _{x%
\protect\varepsilon }^{2}(q)\protect\gamma _{x\protect\varepsilon }(r+q)%
\protect\gamma _{\protect\varepsilon x}(r-p)$}
\end{figure}
\begin{figure}[!h]
\centering \includegraphics{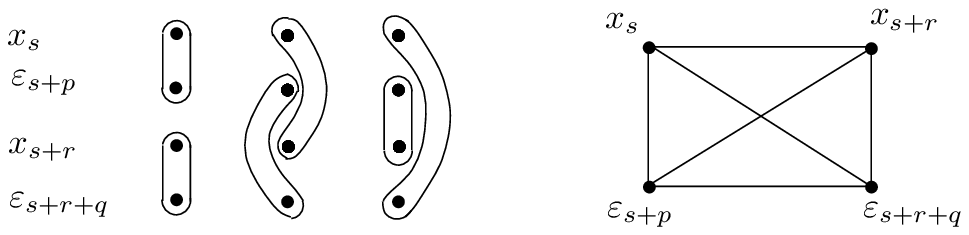} \label{fig:grafo2}
\caption{$\protect\gamma _{x\protect\varepsilon }(r+q)\protect\gamma _{%
\protect\varepsilon x}(r-p)\protect\gamma _{xx}^{2}(r)\protect\gamma _{%
\protect\varepsilon \protect\varepsilon }^{2}(r+q-p)$}
\end{figure}
\begin{figure}[!h]
\centering \includegraphics{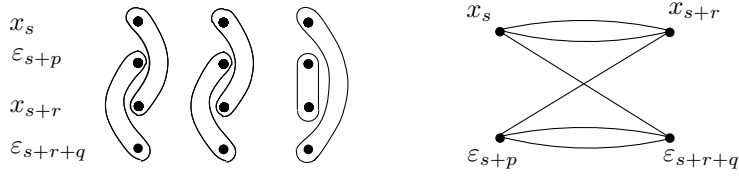} \label{fig:grafo3}
\caption{$\protect\gamma _{x\protect\varepsilon }(p)\protect\gamma _{x%
\protect\varepsilon }(q)\protect\gamma _{x\protect\varepsilon }^{2}(r+q)%
\protect\gamma _{\protect\varepsilon x}^{2}(r-p)\protect\gamma _{xx}(r)%
\protect\gamma _{\protect\varepsilon \protect\varepsilon }(r+q-p)$}
\end{figure}
\begin{figure}[!h]
\centering \includegraphics{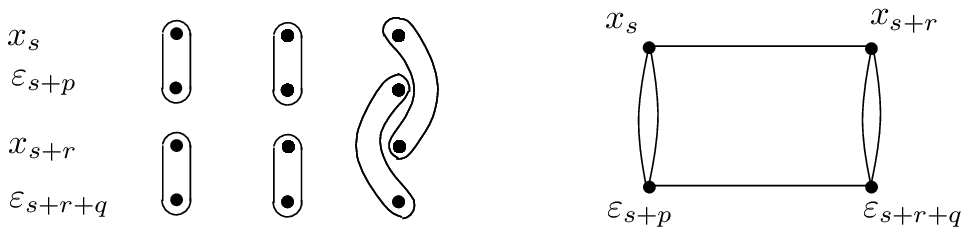} \label{fig:grafo4}
\caption{$\protect\gamma _{x\protect\varepsilon }^{2}(p)\protect\gamma _{x%
\protect\varepsilon }^{2}(q)\protect\gamma _{xx}(r)\protect\gamma _{\protect%
\varepsilon \protect\varepsilon }(r+q-p)$}
\end{figure}
\begin{figure}[!h]
\centering \includegraphics{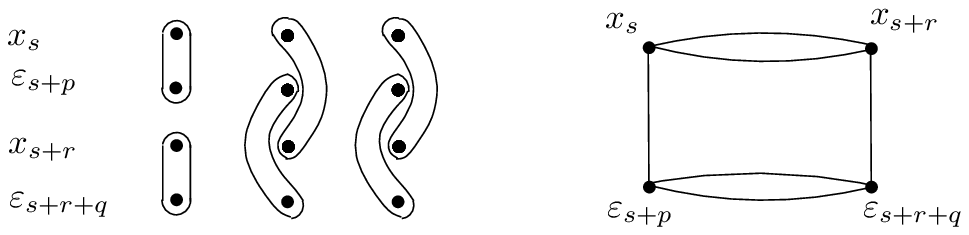} \label{fig:grafo5}
\caption{$\protect\gamma _{x\protect\varepsilon }^{2}(p)\protect\gamma _{x%
\protect\varepsilon }^{2}(q)\protect\gamma _{xx}^{2}(r)\protect\gamma _{%
\protect\varepsilon \protect\varepsilon }^{2}(r+q-p)$}
\end{figure}
\begin{figure}[!h]
\centering \includegraphics{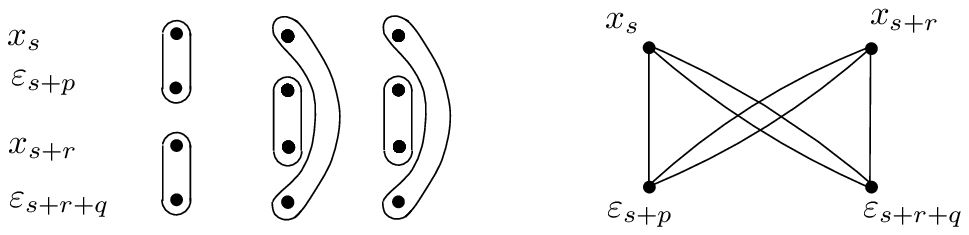} \label{fig:grafo6}
\caption{$\protect\gamma _{x\protect\varepsilon }(p)\protect\gamma _{x%
\protect\varepsilon }(q)\protect\gamma _{x\protect\varepsilon }^{2}(r+q)%
\protect\gamma _{\protect\varepsilon x}^{2}(r-p)$}
\end{figure}
\begin{figure}[!h]
\centering \includegraphics{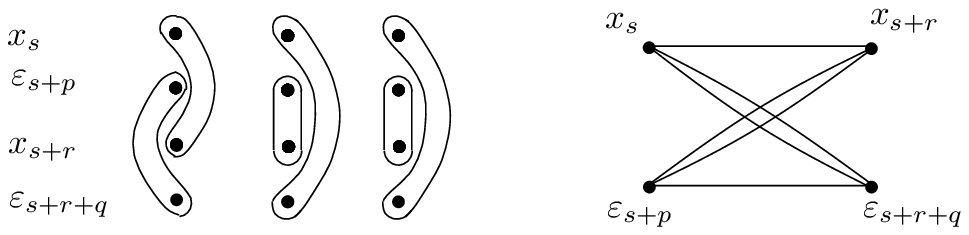} \label{fig:grafo7}
\caption{$\protect\gamma _{xx}(r)\protect\gamma _{\protect\varepsilon 
\protect\varepsilon }(r+q-p)\protect\gamma _{x\protect\varepsilon }^{2}(r+q)%
\protect\gamma _{\protect\varepsilon x}^{2}(r-p)$}
\end{figure}
\bibliographystyle{econometrica}
\bibliography{biblio_NFC}

\ 

Corresponding author:

Domenico Marinucci

Dipartimento di Matematica

Universita' di Roma Tor Vergata

via della Ricerca Scientifica, 1

00133 Roma, Italy

email: marinucc@mat.uniroma2.it

\end{document}